\documentclass[11pt,a4paper,reqno,english]{ijocta}

\usepackage{graphicx}
\usepackage[margin=2cm]{geometry}
\usepackage[font=small,labelfont=bf,tableposition=top]{caption}
\usepackage[font=footnotesize]{subcaption}
\usepackage{multicol}
\usepackage{babel}
\usepackage{amssymb,amsmath,amsfonts}
\usepackage{multirow}
\usepackage{blindtext} 
\usepackage{adjustbox}
\usepackage{hyperref}[breaklinks]
\usepackage{url}
\usepackage{bm}
\usepackage{algorithm}
\usepackage{algpseudocode}

\newtheorem{theorem}{Theorem}

\newtheorem{definition}{Definition}

\newtheorem{lemma}{Lemma}

\newtheorem{property}{Property}
\newtheorem{proposition}{Proposition}

\newcommand{\ceq}{\mathrel{\vcenter{\hbox{:=}}}}

\begin{document}

\title[\small{\textit{Surface Movement Method for Linear Programming}}]{\textbf{Surface Movement Method for Linear Programming}}
\author[\small{\textit{N.A.~Olkhovsky, L.B.~Sokolinsky}}]{\normalsize Nikolay~A.~Olkhovsky, Leonid~B.~Sokolinsky\vspace{0.58cm} 
\\ \small
School of Electronic Engineering and Computer Science, \\ South Ural State University (National Research University), Russia \\
olkhovskiina@susu.ru, leonid.sokolinsky@susu.ru \vspace{0.4cm}}

\begin{abstract}
The article presents a new method of linear programming, called the surface movement method. This method constructs an optimal objective path on the surface of the feasible polytope from the initial boundary point to the point at which the optimal value of the objective function is achieved. The optimality of the path means moving in the direction of maximum increase/decrease in the value of the objective function. A formal description of the algorithm implementing the surface movement method is described. The convergence theorem of this algorithm is proved. The presented method can be effectively implemented using a feed forward deep neural network to determine the optimal direction of movement along the faces of the feasible polytope. To do this, a multidimensional local image of the linear programming problem is constructed at the point of the current approximation. This image is fed to the input of the deep neural network, which returns a vector determining the direction of the optimal objective path on the polytope surface.

\vspace{7pt}
\noindent
\textbf{Keywords: }{Linear programming; Surface movement method; Convergence theorem}

\vspace{2pt}
\noindent
\textbf{AMS Classification:} 90C05
\end{abstract}

\maketitle

\begin{multicols*}{2}

\section{Introduction}

\noindent The rapid development of big data processing technologies in the last decades \cite{Hartung2018:Eng} caused to the emergence of optimization mathematical models in the form of large-scale linear programming (LP) problems~\cite{Sokolinskaya2017a:TR,Veatch2021:Eng}. Of particular interest are the LP problems related to the optimization of non-stationary processes~\cite{Branke2002:Eng}. In a non-stationary LP problem, the constraints and the objective function can change dynamically during computing its solution~\cite{Eremin1979:Rus}. The following optimization problems can be reduced to non-stationary LP problems: choosing the best high-frequency trading strategies~\cite{Brogaard2014:Eng}, optimal navigation and control of aircraft~\cite{Tewari2011:Eng}, dynamic optimization of batch processes~\cite{Srinivasan2003:Eng}, logistics and transportation~\cite{Fleming2021:Eng,Scholl2019:Eng,Meisel2011:Eng}, production planning and control~\cite{Kiran2019:Eng}. Separately, we can mention optimization problems that must be solved in real-time mode~\cite{Mall2007:Eng}. Examples are the following: chemical plant control, energy management, traffic control, multi-point fuel injection system in ICE, autopilot systems, missile guidance system and others. Typically, such applications have time bounds from a few microseconds to several milliseconds.

The simplest approach to solving non-stationary optimization problems is that any changes in the input data is perceived as a separate new problem~\cite{Branke2002:Eng}. This approach may be acceptable when changes occur relatively slowly, and the optimization problem is solved relatively quickly. However, for large-scale non-stationary optimization problems, the solution obtained in this way may be far from optimal due to changes in the input data during computations. In this case, it is necessary to use algorithms that dynamically correct the computational process in accordance with the changing input data. Thus, computations with changed data do not start from scratch, but use information obtained in the past. This approach is applicable to solving optimization problems in real time, provided that the algorithm tracks the movement of the optimal point quickly enough. For large-scale LP problems, the latter requirement makes it urgent to develop scalable methods and parallel algorithms for linear programming.

Up until the present time, one of the most popular ways to solve linear programming problems is a family of algorithms based on the simplex method~\cite{Dantzig1998:Eng}. The simplex method is capable of solving large-scale LP problems by effectively using various types of hyper-sparsity~\cite{Hall2005:Eng}. However, the simplex method has a number of fundamental drawbacks. First, in the worst case, the simplex method must visit all the vertices of the polytope that bound the feasible region, which corresponds to exponential time complexity~\cite{Disser2023:Eng,Hopp2020:Eng,Klee1972:Eng}. Second, the use of the simplex method for solving LP problems with a dimension greater than 50\,000 results in a precision loss~\cite{Bartels1971:Eng} that cannot be corrected by such compute-intensive algorithms as affine scaling or iterative refinement~\cite{Tolla2014:Eng}. Third, the communication structure of algorithms based on the simplex method generally has a limited degree of parallelism, which makes it impossible to efficiently parallelize them on large multiprocessor computing systems with distributed memory~\cite{Hall2010:Eng,Mamalis2015:Eng}. All of the above makes it difficult to use the simplex method to solve non-stationary large-scale LP problems in real-time mode.

Another popular approach to solving large LP problems is a class of algorithms based on the interior-point method~\cite{Zorkaltsev2018:TR}. For the first time, this method was described by Dikin in~\cite{Dantzig2003a:Eng,Dikin1967:TR}. The interior-point method is capable of solving LP problems with millions of variables and constraints~\cite{Gondzio2012:Eng}. The advantage of the interior-point method is that it is self-correcting and is able to provide a high computational accuracy. The main disadvantages of the interior-point method are the following. First, an important subclass of algorithms based on the interior-point method requires some inner point of the feasible region as an initial approximation. Finding such a point can be reduced to solving an supplementary LP problem~\cite{Roos2005:Eng}. Another way to find an inner point is to use Fej\'{e}r approximations~\cite{Sokolinskaya2018:Eng}. Second, the interior-point method does not scale well in large cluster computing systems. There are some special cases when effective parallelization of the interior-point method is possible (see, for example,~\cite{Gondzio2006:Eng}), but in general it is not possible to build an efficient parallel implementation of this method for cluster computing systems. Third, the iterative nature of the interior-point method does not allow predicting the computation time in advance for a certain LP problem. These disadvantages make it difficult to use the interior-point method to solve large-scale LP problems in real-time mode.

Artificial neural networks (ANN)~\cite{Prieto2016:Eng} are a promising new approach to solving optimization problems, which attracts a lot of attention. ANNs is a powerful universal tool that is applicable in almost all problem areas. One of the first applications of ANNs to solving LP problems was the work of Hopfield and Tank~\cite{Tank1986:Eng}. The Hopfield-Tank ANN consists of two fully connected layers and is recurrent. The number of neurons in the first layer is equal to the number of LP problem variables. The number of neurons in the second layer is equal to the number of constraints. The weights and biases of the ANN are completely determined by the parameters of the LP problem. The output is cyclically fed to the input of the ANN. The ANN works until an equilibrium state is reached, when the output becomes equal to the input. This equilibrium state corresponds to a minimum of a special energy function and is a solution to the LP problem. There are a lot of works that extend the Hopfield-Tank approach (see, for example,~\cite{Kennedy1987:Eng,Liu2016:Eng,Malek2005:Eng,Rodriguez-Vazquez1990:Eng,Zak1995:Eng}). The main disadvantage of this approach is that it is impossible to predict the number of ANN work cycles required to achieve the equilibrium state. This makes it impossible to use such recurrent networks to solve large-scale LP problems in real-time mode. For this purpose, feed forward deep neural networks seem more promising. The structure and parameters of such networks, as a rule, do not depend on the particular input parameters of the problem. The solution is obtained in one pass with a fixed network operation time, which makes it possible to use them to solve problems in real-time mode. One of the important classes of ANNs are convolutional neural networks~\cite{LeCun2015:Eng}. This class is of special interest for image processing. In recent paper~\cite{Olkhovsky2022:TR}, a new method for constructing images of multidimensional LP problems was proposed, which makes it possible to use feed forward neural networks, including convolutional ones, to solve them. It should be noted that deep neural networks require training on a large number of labeled datasets, which can be efficiently performed on GPUs~\cite{Raina2009:Eng}. In paper~\cite{Sokolinsky2023:TR}, the apex method is proposed for solving LP problems, which makes it possible to construct a path on the surface of the feasible polytope in the direction of maximum increase/decrease in the value of the objective function, leading to the optimum point. The apex method belongs to the class of iterative projection-type methods, which are characterized by a low linear convergence rate, which makes them unacceptable for real-time mode.

This article presents a new method for solving LP problems, called the surface movement method. This method is intended for using feed forward ANNs, including convolutional neural networks. The rest of the paper is organized as follows. Section~\ref{sec:theoretical_basis} presents the theoretical background on which the surface movement method is based. Section~\ref{sec:surface-movement-method-description} contains a description of the surface movement method and a proof of the convergence theorem. In Section~\ref{sec:discussion}, we discuss the strengths and weaknesses of the proposed method, as well as reveal the ways of its practical implementation based on the synthesis of supercomputer and neural network technologies. Section~\ref{sec:conclusion} summarizes the results and provides further research directions. A summary of the main symbols used in the paper is presented in Appendix.

\section{Theoretical Background}\label{sec:theoretical_basis}

\noindent This section presents a theoretical background on which the surface movement method is based. We consider a LP problem in the following form:
\begin{equation}\label{eq:LP_problem}
	\bm{\bar{x}} = \arg \mathop {{\max}}\limits_{\bm{x}\in\mathbb{R}^n} \left\{ {\left\langle \bm{c},\bm{x} \right\rangle \left| A\bm{x} \leqslant \bm{b}\right.}  \right\},	
\end{equation}
where $\bm{c}\in \mathbb{R}^n$, $\bm{b} \in \mathbb{R}^m$, $A \in \mathbb{R}^{m\times n}$, $m>1$, $\bm{c}\ne \mathbf{0}$. Here, $\left\langle \cdot,\cdot \right\rangle$ stands for the dot product of two vectors. We assume that the constraint $\bm{x} \geqslant \mathbf{0}$ is also included in the matrix inequality $A\bm{x} \leqslant \bm{b}$ in the form of $-\bm{x} \leqslant \mathbf{0}$. Denote by $\mathcal P$ the set of indexes numbering the rows of the matrix $A$:
\begin{equation*}
	\mathcal{P}=\left\lbrace 1,\cdots,m\right\rbrace.
\end{equation*}
The linear objective function of the problem~\eqref{eq:LP_problem} has the form
\begin{equation*}
	f(\bm{x})=\left\langle \bm{c},\bm{x}\right\rangle.
\end{equation*}
In this case, the vector $\bm{c}$ is the gradient of the objective function $f(\bm{x})$.

Let $\bm{a}_i\in\mathbb{R}^n$ denote a vector representing the $i$th row of the matrix~$A$. We assume that $\bm{a}_i\ne\mathbf{0}$ for all $i\in\mathcal P$. Denote by $\hat H_i$ a closed half-space defined by the inequality $\left\langle\bm{a}_i,\bm{x}\right\rangle\leqslant b_i$, and by $H_i$ --- the hyperplane bounding it:
\begin{equation}\label{eq:half-space}
	\hat H_i=\left\lbrace \bm{x}\in\mathbb{R}^n \middle|\left\langle  \bm{a}_i,\bm{x} \right\rangle \leqslant b_i \right\rbrace;
\end{equation}
\begin{equation}\label{eq:H_i}
	H_i=\left\lbrace \bm{x}\in\mathbb{R}^n \middle|\left\langle  \bm{a}_i,\bm{x} \right\rangle = b_i \right\rbrace.
\end{equation}
Let us define a feasible polytope
\begin{equation}\label{eq:polytope_M}
	M = \bigcap\limits_{i \in \mathcal{P}} \hat H_i,
\end{equation}
representing the feasible region of LP~Problem~\eqref{eq:LP_problem}. Note that $M$, in this case, is a closed convex set. We assume that $M$ is bounded, and $M\ne\emptyset$, i.e., LP~Problem~\eqref{eq:LP_problem} has a solution. Denote by $\Gamma(M)$ the set of boundary points of the polytope $M$\footnote{If $M\subset\mathbb{R}^n$, then a point $\bm{x} \in \mathbb{R}^n$ is a boundary point of $M$ if every neighborhood of $\bm{x}$ contains at least one point in $M$ and at least one point not in $M$.}.

\begin{proposition}\label{prp:Gamma(M)}
	Let $M$ be the feasible polytope of LP~Problem~\eqref{eq:LP_problem}, defined by equation~\eqref{eq:polytope_M}. Let $V_\epsilon({\bm{u}})$ denotes the $\epsilon$-neighborhood of the point ${\bm{u}}$ in $\mathbb{R}^n$. Then, 
	\begin{eqnarray*}
		&& \forall {\bm{u}} \in \Gamma(M)~\exists \epsilon > 0 :  \\
		&& \forall \bm{w} \in V_\epsilon({\bm{u}}) \cap \Gamma(M) \exists i' \in \mathcal{P}: {\bm{u}},\bm{w} \in H_{i'},
	\end{eqnarray*}
	i.e., for any point ${\bm{u}}\in \Gamma(M)$, there exists $\epsilon>0$ such that for each boundary point $\bm{w}$ belonging to the $\epsilon$-neighborhood of the point ${\bm{u}}$, there exists at least one $i'\in \mathcal{P}$ such that ${\bm{u}},\bm{w}\in H_{i'}$.
\end{proposition}
\begin{proof}
	Fix an arbitrary point ${\bm{u}} \in \Gamma(M)$. Define
	\begin{equation}\label{eq:P_u}
		\mathcal{P}_{\bm{u}} = \left\lbrace \left. i \in \mathcal{P} \right| {\bm{u}} \in H_i \right\rbrace.
	\end{equation}
	In other words, $\mathcal{P}_{\bm{u}}$ is the set of indices of all hyperplanes $H_i$ to which the point~${\bm{u}}$ belongs. Denote
	\begin{equation*}
		\mathcal{P}_{\backslash {\bm{u}}} = \mathcal{P} \backslash \mathcal{P}_{\bm{u}},
	\end{equation*}
	i.e., $\mathcal{P}_{\backslash {\bm{u}}}$ is the set of indices of all hyperplanes $H_i$ to which the point~${\bm{u}}$ does not belong. Define
	\begin{equation*}
		\delta = \min \left\lbrace \left.  dist({\bm{u}},H_i) \right| i \in \mathcal{P}_{\backslash {\bm{u}}} \right\rbrace,
	\end{equation*}
	where $dist({\bm{u}},H_i)$ stands for the Euclidean distance from the point ${\bm{u}}$ to the hyperplane $H_i$\footnote{In this case, $dist({\bm{u}},H_i)=\frac{\left\langle \bm{a}_i,{\bm{u}}\right\rangle - b_i}{\left\|\bm{a}_i\right\|}$.}. By definition, 
	\begin{equation*}
		\delta > 0.
	\end{equation*}
	Take $\epsilon$ satisfying the condition
	\begin{equation*}
		0<\epsilon<\delta.
	\end{equation*}
	Then for any $\bm{w} \in V_\epsilon({\bm{u}}) \cap \Gamma(M)$, the following condition holds:
	\begin{equation*}
		\forall i \in \mathcal{P}_{\backslash {\bm{u}}} : \bm{w} \notin H_i.
	\end{equation*}
	Since $\bm{w}$ is a boundary point, it follows that there exists an $i' \in \mathcal{P}_{\bm{u}}$ such that
	\begin{equation*}
		\bm{w} \in H_{i'}.
	\end{equation*}
	By virtue of~\eqref{eq:P_u}, the following condition also holds:
	\begin{equation*}
		{\bm{u}} \in H_{i'}.
	\end{equation*}
\end{proof}

\begin{definition}\label{def:c-projection_onto_H_i}
	The objective projection of a point $\bm{z}\in\mathbb{R}^n$ onto the hyperplane $H_i$ is a point $\bm{\gamma}_i(\bm{z})\in\mathbb{R}^n \cup \left\lbrace \infty \right\rbrace$ defined by the equation
	\begin{equation}\label{eq:gamma_i0}
		\bm{\gamma}_i(\bm{z}) = \left\{ 
		\begin{array}{l}
			L(\bm{z})\cap H_i,\,\text{if}~\left\langle \bm{a}_i,\bm{c}\right\rangle \ne 0;\\
			\infty,~~~~~~~~~~~\text{if}~\left\langle \bm{a}_i,\bm{c}\right\rangle = 0,
		\end{array} 
		\right.
	\end{equation}
	where $L(\bm{z})$ is the line passing through the point $\bm{z}$ parallel to the vector~$\bm{c}$:
	\begin{equation}\label{eq:L(z)}
		L(\bm{z})=\left\lbrace \left. {\bm{y}}\in\mathbb{R}^n \right| {\bm{y}}=\bm{z}+\lambda \bm{c}, \lambda\in\mathbb{R}\right\rbrace.
	\end{equation}
	In other words, if the vector $\bm{c}$ is not parallel to the hyperplane $H_i$, then the objective projection of the point $\bm{z}$ onto the hyperplane $H_i$ is the intersection point of this hyperplane with the line passing through the point $\bm{z}$ parallel to the vector~$\bm{c}$. In the case when the vector $\bm{c}$ is parallel to the hyperplane $H_i$, the objective projection is the point at infinity.
\end{definition}

The following proposition provides an equation for calculating the objective projection $\bm{\gamma}_i(\bm{z})$ of point $\bm{z}$ onto the hyperplane $H_i$.

\begin{proposition}\label{prp:c-projection}
	Let $\left\langle \bm{a}_i,\bm{c}\right\rangle \ne 0$, i.e., the vector~$\bm{c}$ is not parallel to the hyperplane $H_i$. Then
	\begin{equation}\label{eq:gamma_i(z)}
		\bm{\gamma}_i(\bm{z}) =  \bm{z} - \frac{\left\langle \bm{a}_i,\bm{z}\right\rangle - b_i}{\left\langle \bm{a}_i,\bm{c}\right\rangle}\bm{c}.
	\end{equation}
\end{proposition}
\begin{proof}
	In accordance with~\eqref{eq:gamma_i0} and~\eqref{eq:L(z)}, the following equation holds for some $\lambda\in\mathbb{R}$:
	\begin{equation}\label{eq:gamma=z+lambda*c}
		\bm{\gamma}_i(\bm{z}) = \bm{z}+\lambda \bm{c}
	\end{equation}
	On the other hand, in accordance with~\eqref{eq:H_i}, the following equation holds:
	\begin{equation}\label{eq:<a_i,gamma>=b_i}
		\left\langle  \bm{a}_i,\bm{\gamma}_i(\bm{z}) \right\rangle = b_i.
	\end{equation}
	Substituting the right-hand side of equation~\eqref{eq:gamma=z+lambda*c} instead of~$\bm{\gamma}_i(z)$ into equation~\eqref{eq:<a_i,gamma>=b_i}, we obtain
	\begin{equation*}
		\left\langle  \bm{a}_i,\bm{z}+\lambda \bm{c} \right\rangle = b_i.
	\end{equation*}
	It follows that
	\begin{equation}\label{eq:lambda}
		\lambda=-\frac{\left\langle  \bm{a}_i,\bm{z} \right\rangle -b_i}{\left\langle  \bm{a}_i,\bm{c} \right\rangle}.
	\end{equation}
	Substituting the right-hand side of equation~\eqref{eq:lambda} instead of~$\lambda$ into equation~\eqref{eq:gamma=z+lambda*c}, we obtain
	\begin{equation*}
		\bm{\gamma}_i(\bm{z}) =  \bm{z} - \frac{\left\langle \bm{a}_i,\bm{z}\right\rangle - b_i}{\left\langle \bm{a}_i,\bm{c}\right\rangle}\bm{c}.
	\end{equation*}
\end{proof}

\begin{definition}\label{def:bias}
	The objective bias of the point $\bm{z}\in\mathbb{R}^n$ relative to the hyperplane $H_i$ is the scalar quantity~$\beta_i(\bm{z})$, calculated by the equation
	\begin{equation}\label{eq:beta_i}
		\beta_i(\bm{z})=-\frac{\left\langle \bm{a}_i,\bm{z}\right\rangle - b_i}{\left\langle \bm{a}_i,\bm{c}\right\rangle}\left\| \bm{c} \right\|.
	\end{equation}
\end{definition}

\noindent For brevity's sake, everywhere below, we will use the term ``bias'', meaning by this ``objective bias''. Denote
\begin{equation}\label{eq:e_c}
	\bm{e_c} = \frac{\bm{c}}{\left\|\bm{c}\right\| }.
\end{equation}
Then equation~\eqref{eq:gamma_i(z)} can be rewritten as follows:
\begin{equation}\label{eq:gamma_i}
	\bm{\gamma}_i(\bm{z}) =  \bm{z} + \beta_i(\bm{z})\bm{e_c},
\end{equation}
which is equivalent to
\begin{equation*}
	\beta_i(\bm{z})\bm{e_c} =  \bm{\gamma}_i(\bm{z}) - \bm{z}.
\end{equation*}
Taking into account~\eqref{eq:e_c}, it follows
\begin{equation}\label{eq:|beta_i(z)|}
	\left| \beta_i(\bm{z})\right|  =  \left\| \bm{\gamma}_i(\bm{z}) - \bm{z}\right\|.
\end{equation}
Thus, $\left| \beta_i(\bm{z})\right|$ is the distance from the point $\bm{z}$ to its objective projection onto the hyperplane $H_i$.

\begin{definition}\label{def:H_c} 
	The objective hyperplane $H_c(\bm{z})$ passing through the point $\bm{z}$ is the hyperplane defined by the following equation:
	\begin{equation}\label{eq:H_c} 
		H_c(\bm{z})=\left\lbrace \bm{x}\in\mathbb{R}^n \middle|\left\langle  \bm{c},\bm{x} \right\rangle = \left\langle  \bm{c},\bm{z} \right\rangle \right\rbrace.
	\end{equation}
\end{definition}

The following proposition holds.

\begin{proposition}\label{prp:metrics}
	Fix an arbitrary point ${\bm{z}\in\mathbb{R}^n}$. Then, for any points \[{\bm{z}',\bm{z}'' \in H_c(\bm{z})}, \bm{z}'\ne \bm{z}'',\] the following statement is true for all $i\in\mathcal{P}$:
	\begin{equation*}
		\left\langle \bm{c},\bm{\gamma}_i(\bm{z}')\right\rangle < \left\langle \bm{c},\bm{\gamma}_i(\bm{z}'')\right\rangle \Leftrightarrow \beta_i(\bm{z}') < \beta_i(\bm{z}'').
	\end{equation*}
\end{proposition}
\begin{proof}
	Since $\bm{z}',\bm{z}'' \in H_c(\bm{z})$, and $\bm{c}$ is normal to the hyperplane $H_c(\bm{z})$, the following two equations hold:
	\begin{equation*}
		\left\langle \bm{c}, \bm{z}'-\bm{z}\right\rangle = 0;
	\end{equation*}
	\begin{equation*}
		\left\langle \bm{c}, \bm{z}''-\bm{z}\right\rangle = 0.
	\end{equation*}
	Hence,
	\begin{equation}\label{eq:<c,z''>=<c,z>=<c,z'>}
		\left\langle \bm{c}, \bm{z}''\right\rangle = \left\langle \bm{c}, \bm{z}\right\rangle = \left\langle \bm{c}, \bm{z}'\right\rangle.
	\end{equation}
	Sequentially using~\eqref{eq:gamma_i}, \eqref{eq:<c,z''>=<c,z>=<c,z'>}, and \eqref{eq:e_c}, we obtain the following chain of equivalent inequalities:
	\begin{eqnarray*}
		&& \left\langle \bm{c},\bm{\gamma}_i(\bm{z}')\right\rangle < \left\langle \bm{c},\bm{\gamma}_i(\bm{z}'')\right\rangle \\ 
		&& \Leftrightarrow \left\langle \bm{c},\bm{z}' + \beta_i(\bm{z}')\bm{e_c}\right\rangle < \left\langle \bm{c},\bm{z}'' + \beta_i(\bm{z}'')\bm{e_c}\right\rangle\\
		&& \Leftrightarrow \left\langle \bm{c},\beta_i(\bm{z}')\bm{e_c}\right\rangle < \left\langle \bm{c},\beta_i(\bm{z}'')\bm{e_c}\right\rangle\\
		&& \Leftrightarrow \left\langle \bm{c},\beta_i(\bm{z}')\bm{c}/\left\| \bm{c}\right\|\right\rangle < \left\langle \bm{c},\beta_i(\bm{z}'')\bm{c}/\left\| \bm{c}\right\| \right\rangle \\
		&& \Leftrightarrow \frac{\beta_i(\bm{z}')}{\left\| \bm{c}\right\|}\left\langle \bm{c},\bm{c}\right\rangle < \frac{\beta_i(\bm{z}'')}{\left\| \bm{c}\right\|}\left\langle \bm{c},\bm{c}\right\rangle \\
		&& \Leftrightarrow \beta_i(\bm{z}') < \beta_i(\bm{z}'').
	\end{eqnarray*}
\end{proof}

The following definition of a recessive half-space is adopted by us from the paper~\cite{Sokolinsky2023:TR}.

\begin{definition}\label{def:recessive_halfspace} 
	The half-space $\hat H_i$ is called recessive if
	\begin{equation}\label{eq:recessive_half-space}
		\forall \bm{x}\in {H_i},\forall \lambda>0:\bm{x}+\lambda \bm{c}\notin \hat H_i.
	\end{equation}
	The geometric interpretation of this definition is as follows: the ray parallel to the vector~$\bm{c}$ coming from any point of the hyperplane bounding the recessive half-space has no points in common with this half-space, except for the starting one.
\end{definition}

The following condition is both necessary and sufficient for the half-space~$\hat H_i$ to be recessive~\cite{Sokolinsky2023:TR}:
\begin{equation}\label{eq:recessive_half-space_condition}
	\left\langle \bm{a}_i,\bm{c}\right\rangle>0.
\end{equation}

A recessive half-space has the following properties.
\begin{property}\label{pty:recessive_property1}
	Let the half-space $\hat{H}_i$ be recessive. Then, any line parallel to the vector $\bm{c}$ intersects the hyperplane $H_i$ at a single point.
\end{property}
\noindent This property directly follows from the fact that the hyperplane $H_i$ bounding the recessive half-space $\hat{H}_i$ cannot be parallel to the vector $\bm{c}$ by Definition~\ref{def:recessive_halfspace}.

\begin{property}\label{pty:recessive_property2}
	Let the half-space $\hat{H}_i$ be recessive. Then, 
	\begin{equation}\label{eq:property2}
		\bm{x}\in\hat{H}_i \Leftrightarrow \beta_i(\bm{x})\geqslant 0.
	\end{equation}
\end{property}
\begin{proof} 
	First, assume that $\bm{x}\in\hat{H}_i$. Then, according to~\eqref{eq:half-space}, the following inequality holds:
	\begin{equation*}
		\left\langle \bm{a}_i,\bm{x}\right\rangle - b_i\leqslant 0.
	\end{equation*}
	By virtue of~\eqref{eq:beta_i}, we have
	\begin{equation}\label{eq:beta_i(x)}
		\beta_i(\bm{x})=-\frac{\left\langle \bm{a}_i,\bm{x}\right\rangle - b_i}{\left\langle \bm{a}_i,\bm{c}\right\rangle}\left\| \bm{c} \right\|.
	\end{equation}
	Taking into account~\eqref{eq:recessive_half-space_condition}, we obtain
	\begin{equation*}
		\bm{x}\in\hat{H}_i \Rightarrow \beta_i(\bm{x})\geqslant 0.
	\end{equation*}
	
	Now, suppose that $\beta_i(\bm{x})\geqslant 0$. According to~\eqref{eq:beta_i}, this means that
	\begin{equation*}
		\frac{\left\langle \bm{a}_i,\bm{x}\right\rangle - b_i}{\left\langle \bm{a}_i,\bm{c}\right\rangle}\left\| \bm{c} \right\| \leqslant 0.
	\end{equation*}
	Given~\eqref{eq:recessive_half-space_condition}, we obtain from here
	\begin{equation*}
		\left\langle \bm{a}_i,\bm{x}\right\rangle - b_i \leqslant 0.
	\end{equation*}
	According to~\eqref{eq:half-space}, it follows that
	\begin{equation*}
		\bm{x}\in\hat{H}_i.
	\end{equation*}
	Thus,
	\begin{equation*}
		\beta_i(\bm{x})\geqslant 0 \Rightarrow \bm{x}\in\hat{H}_i.
	\end{equation*}
\end{proof}

Define
\begin{equation}\label{eq:I}
	\mathcal{I}=\left\lbrace i\in\mathcal{P} \left| \left\langle \bm{a}_i,\bm{c}\right\rangle >0 \right. \right\rbrace,
\end{equation}
i.e., $\mathcal{I}$ represents a set of indices for which the half-space $\hat H_i$ is recessive. Since the feasible polytope $M$ is a bounded set, we have
\begin{equation}\label{Ic_ne_eptyset}
	\mathcal{I}\ne\emptyset.
\end{equation}
Define
\begin{equation}\label{eq:polytope_hat_M}
	\hat{M} = \bigcap\limits_{i \in \mathcal{I}} \hat H_i.
\end{equation}
Obviously, $\hat{M}$ is a convex, closed, unbounded polytope. Let us call it a recessive polytope. By~\eqref{eq:polytope_M} and~\eqref{eq:I}, it follows that
\begin{equation}\label{eq:M_subset_hat_M}
	M \subset \hat{M}.
\end{equation}
Let $\Gamma(\hat M)$ denote the set of boundary points of the recessive polytope $\hat{M}$. According to Proposition~3 in~\cite{Sokolinsky2023:TR}, we have
\begin{equation*}
	\bar{\bm{x}} \in \Gamma(\hat M),
\end{equation*}
i.e., the solution of LP~Problem~\eqref{eq:LP_problem} lies on the boundary of the recessive polytope $\hat{M}$.

\begin{proposition}\label{prp:Gamma(hat_M)}
	Let $\hat{M}$ be the recessive polytope defined by equation~\eqref{eq:polytope_hat_M}. Then for any point ${\bm{u}} \in \Gamma(\hat{M})$, there exists $\epsilon > 0$ such that for any boundary point~$\bm{w}$ belonging to the $\epsilon$-neighborhood $V_\epsilon({\bm{u}})$ of the point~${\bm{u}}$, there is $i' \in \mathcal{I}$ for which ${\bm{u}},\bm{w} \in H_{i'}$ is valid, i.e.,
	\begin{eqnarray*}
		&& \forall {\bm{u}} \in \Gamma(\hat{M})~\exists \epsilon > 0 : \\
		&& \forall \bm{w} \in V_\epsilon({\bm{u}}) \cap \Gamma(\hat{M})~\exists i' \in \mathcal{I}: {\bm{u}},\bm{w} \in H_{i'}.
	\end{eqnarray*}
\end{proposition}
\begin{proof} The proof is similar to the proof of proposition~\ref{prp:Gamma(M)}.
\end{proof}	

\begin{definition}\label{def:objective_projection_onto_hat_G} 
	The objective projection of the point $\bm{z}\in\mathbb{R}^n$ onto the boundary $\Gamma(\hat M)$ of the recessive polytope $\hat{M}$ is the point $\hat{\bm{\gamma}}(\bm{z})$ calculated by the equation
	\begin{equation*}
		\hat{\bm{\gamma}}(\bm{z})=L(\bm{z}) \cap \Gamma(\hat M),
	\end{equation*}
	where $L(\bm{z})$ --- the line passing through the point $\bm{z}$ parallel to the vector $\bm{c}$:
	\begin{equation*}
		L(\bm{z})=\left\lbrace \left. {\bm{y}}\in\mathbb{R}^n \right| {\bm{y}}=\bm{z}+\lambda \bm{c}, \lambda\in\mathbb{R}\right\rbrace.
	\end{equation*}
	The scalar quantity $\hat\beta(\bm{z})\in\mathbb{R}$ satisfying the equation
	\begin{equation}\label{eq:hat_gamma(z)=z+beta(z)c}
		\hat{\bm{\gamma}}(\bm{z})=\bm{z}+\hat\beta(\bm{z})\bm{c}
	\end{equation}
	will be called the bias of the point $\bm{z}$ relative to the boundary of the recessive polytope $\hat M$.
\end{definition}

\noindent Note that the correctness of this definition is based on property~\ref{pty:recessive_property1}. The following proposition provides an equation for calculating the objective projection onto the boundary of the recessive polytope. Fig.~\ref{fig:Proposition5} illustrates the proof of this proposition.

\noindent
\begin{minipage}{\linewidth}
	\centering
	\includegraphics[width=6 cm]{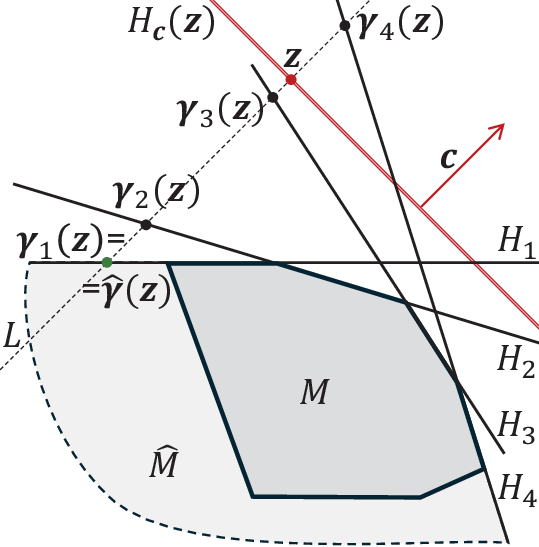}
	\captionof{figure}{Illustration to proof of Proposition \ref{prp:hat_gamma(z)=gamma_i'(z)}.}\label{fig:Proposition5}
\end{minipage}
\medskip

\begin{proposition}\label{prp:hat_gamma(z)=gamma_i'(z)}
	Let an arbitrary point $\bm{z}\in\mathbb{R}^n$ be given. Put
	\begin{equation}\label{eq:arg_min_beta_i(z)}
		i' = \arg\min \left\{ \beta_i(\bm{z})   \left| i\in\mathcal{I} \right.\right\}.	
	\end{equation}
	Then
	\begin{equation}\label{eq:hat_gamma(z)=gamma_i'(z)}
		\hat{\bm{\gamma}}(\bm{z})=\bm{\gamma}_{i'}(\bm{z}).
	\end{equation}
	In other words, the objective projection of the point $\bm{z}$ onto the boundary of the recessive polytope~$\hat M$ coincides with the projection of this point onto the hyperplane $H_{i'}$, which has the minimum bias relative to $\bm{z}$.
\end{proposition}
\begin{proof}
	Fix an arbitrary point $\bm{z}\in\mathbb{R}^n$. According to Definition~\ref{def:objective_projection_onto_hat_G}, let us construct a line parallel to the vector~$\bm{c}$, which passes through the point~$\bm{z}$:
	\begin{equation*}
		L=\left\lbrace \left. {\bm{y}}\in\mathbb{R}^n \right| {\bm{y}}=\bm{z}+\lambda \bm{c}, \lambda\in\mathbb{R}\right\rbrace.
	\end{equation*}
	The following equation holds:
	\begin{equation*}
		\hat{\bm{\gamma}}(\bm{z}) = L\cap \Gamma(\hat{M}).
	\end{equation*}
	In accordance with property \ref{pty:recessive_property1}, for any $i\in\mathcal{I}$, the line $L$ intersects the hyperplane $H_i$ at only one point, which we denote as ${\bm{y}}_i$:
	\begin{equation*}
		L\cap H_i = \left\lbrace {\bm{y}}_i \right\rbrace.
	\end{equation*}
	Define
	\begin{equation*}
		Y=\bigcup\limits_{i\in\mathcal{I}}{}\left\lbrace {\bm{y}}_i \right\rbrace,
	\end{equation*}
	i.e., $Y$ is the set of points at which the line $L$ intersects the boundaries of recessive half-spaces. According to Definition~\ref{def:c-projection_onto_H_i}, the following conditions hold:
	\begin{equation}\label{eq:y_i}
		\forall i \in \mathcal{I} ~:~ \bm{\gamma}_i(\bm{z})={\bm{y}}_i,
	\end{equation}
	and
	\begin{equation}\label{eq:y_i_in_H_i}
		\forall i \in \mathcal{I} ~:~ {\bm{y}}_i\in H_i.
	\end{equation}
	By virtue of \eqref{eq:polytope_hat_M}, the following is also true:
	\begin{equation*}
		\hat{\bm{\gamma}}(\bm{z}) \in Y.
	\end{equation*}
	This means that there is $i' \in \mathcal{I}$ such that
	\begin{equation}\label{eq:gamma^(z)=gamma_i'(z)}
		\hat{\bm{\gamma}}(\bm{z}) = \bm{\gamma}_{i'}(\bm{z}).
	\end{equation}
	Let us show that
	\begin{equation*}
		i' \in \mathrm{Arg}\min \left\{ \beta_i(\bm{z})   \left| i\in\mathcal{I} \right.\right\}.
	\end{equation*}
	Assume the opposite, namely that
	\begin{equation*}
		\beta_{i'}(\bm{z}) > \min \left\{ \beta_i(\bm{z})   \left| i\in\mathcal{I} \right.\right\}.
	\end{equation*}
	Then there exists an $i''\in\mathcal{I}$ such that
	\begin{equation}\label{eq:beta_i''<beta_i'}
		\beta_{i''}(\bm{z})<\beta_{i'}(\bm{z})
	\end{equation}
	By \eqref{eq:gamma_i} and \eqref{eq:y_i}
	\begin{equation*}
		{\bm{y}}_{i'} =\bm{z}+\beta_{i'}(\bm{z})\bm{e_c};
	\end{equation*}
	\begin{equation*}
		{\bm{y}}_{i''}=\bm{z}+\beta_{i''}(\bm{z})\bm{e_c}.
	\end{equation*}
	Hence
	\begin{equation*}
		{\bm{y}}_{i'} = {\bm{y}}_{i''} + \left( \beta_{i'}(\bm{z}) - \beta_{i''}(\bm{z})\right) \bm{e_c},
	\end{equation*}
	which, by virtue of \eqref{eq:beta_i''<beta_i'} and \eqref{eq:e_c} is equivalent to
	\begin{equation}\label{eq:y''}
		{\bm{y}}_{i'} = {\bm{y}}_{i''}+\frac{\left| \beta_{i''}(\bm{z})-\beta_{i'}(\bm{z})\right| }{\left\|\bm{c}\right\| }\bm{c}.
	\end{equation}
	In accordance with \eqref{eq:recessive_half-space} and \eqref{eq:y_i_in_H_i}, it follows from \eqref{eq:y''} that
	\begin{equation*}
		{\bm{y}}_{i'} \notin \hat{H}_{i''}.
	\end{equation*}
	This means that 
	\begin{equation*}
		{\bm{y}}_{i'} \notin \hat{M}.
	\end{equation*}
	Taking into account \eqref{eq:y_i} and \eqref{eq:gamma^(z)=gamma_i'(z)}, it follows that
	\begin{equation*}
		\hat{\bm{\gamma}}(\bm{z}) \notin \hat{M}.
	\end{equation*}
	We have thus reached a contradiction to Definition \ref{def:objective_projection_onto_hat_G}.
\end{proof}

The following proposition provides an equation for calculating the bias of a point relative to the boundary of a recessive polytope.

\begin{proposition}
	Let an arbitrary point $\bm{z}\in\mathbb{R}^n$ be given. Then
	\begin{equation}\label{eq:beta(z)}
		\hat\beta(\bm{z})=\frac{\left\langle  \bm{c},\hat{\bm{\gamma}}(\bm{z})-\bm{z}\right\rangle}{{\left\|  \bm{c} \right\|}^2}.
	\end{equation}
\end{proposition}
\begin{proof}
	In accordance with \eqref{eq:hat_gamma(z)=z+beta(z)c}, the points $\hat{\bm{\gamma}}(\bm{z})$ and $\bm{z}$ are on the normal to the hyperplane $H_c(\bm{z})$. Therefore, the point $\bm{z}\in H_c(\bm{z})$ is an orthogonal projection of the point $\hat{\bm{\gamma}}(\bm{z})$ onto the hyperplane $H_c(\bm{z})$, and taking into account~\eqref{eq:H_c}, it can be calculated by the following known equation:
	\begin{equation}\label{eq:projection_onto_H_c}	
		\bm{z} =\hat{\bm{\gamma}}(\bm{z}) - \frac{\left\langle  \bm{c},\hat{\bm{\gamma}}(\bm{z})-\bm{z}\right\rangle}{{\left\|  \bm{c} \right\|}^2}  \bm{c}.
	\end{equation}
	We can rewrite this as $	\hat{\bm{\gamma}}(\bm{z}) = \bm{z} + \frac{\left\langle  \bm{c},\hat{\bm{\gamma}}(\bm{z})-\bm{z}\right\rangle}{{\left\|  \bm{c} \right\|}^2}\bm{c}$.
	Comparing this with \eqref{eq:hat_gamma(z)=z+beta(z)c}, we obtain
	\begin{equation*}
		\hat\beta(\bm{z})=\frac{\left\langle  \bm{c},\hat{\bm{\gamma}}(\bm{z})-\bm{z}\right\rangle}{{\left\|  \bm{c} \right\|}^2}.
	\end{equation*}\end{proof}

\section{Surface Movement Method}\label{sec:surface-movement-method-description}

The surface movement method constructs, on the surface of the feasible polytope, a path from an arbitrary boundary point ${\bm{u}}^{(0)}\in M \cap\Gamma(\hat{M})$ to a point $\bar{\bm{x}}$, which is a solution to LP~Problem~\eqref{eq:LP_problem}. Moving along the surface of the recessive polytope is performed in the direction of the greatest increase in the value of the objective function. The path constructed as a result of such a movement will be called the optimal objective path. The implementation of the surface movement method is presented in the form of Algorithm~\ref{alg:SMM}. Let us give brief comments on this implementation. Step 1 reads the initial approximation of ${\bm{u}}^{(0)}$. It can be an arbitrary boundary point of a recessive polytope $\hat M$ satisfying the condition
\begin{equation*}
	{\bm{u}}^{(0)}\in M \cap\Gamma(\hat{M}),
\end{equation*}
which is checked in Step~2. 
Step~3 assigns the iteration counter~$k$ the value~0. 
Step~4 sets the value of the parameter~$r$. 
Step~5 builds an $n$-dimensional disk~$D^{(0)}$, which is the intersection of the objective hyperplane~$H_c\left({\bm{u}}^{(0)}\right)$ passing through the point~${\bm{u}}^{(0)}$, and the $n$-dimensional ball $V_r\left({\bm{u}}^{(0)}\right)$ of a small radius~$r$ with the center at the point~${\bm{u}}^{(0)}$. 
Step~6 calculates the point ${\bm{v}}^{(0)}\in D^{(0)}$, having the maximum bias relative to the boundary of the recessive polytope~$\hat{M}$.
\end{multicols}
\begin{minipage}{0.95\textwidth}
\begin{algorithm}[H]\caption{Surface movement method}\label{alg:SMM}
	\begin{algorithmic}[1]
		\Require \begin{flushleft} $\hat H_i=\left\lbrace \bm{x}\in\mathbb{R}^n \middle|\left\langle  \bm{a}_i,\bm{x} \right\rangle \leqslant b_i \right\rbrace;\; \hat M=\bigcap\limits_{i \in \mathcal{I}} {\hat H_i};\; \forall \bm{z} \in \mathbb{R}^n : H_c(\bm{z})=\left\lbrace \bm{x}\in\mathbb{R}^n \middle|\left\langle  \bm{c},\bm{x}-\bm{z} \right\rangle = 0 \right\rbrace$ \end{flushleft}
		\State \textbf{input} ${\bm{u}}^{(0)}$
		\State \textbf{assert} ${\bm{u}}^{(0)}\in M \cap \Gamma(\hat{M})$
		\State $k\ceq 0$
		\State $r\ceq 0.1$
		\State $D^{(0)}\ceq H_c({\bm{u}}^{(0)}) \cap V_r ({\bm{u}}^{(0)})$ \Comment $V_r({\bm{u}}^{(0)})$ is $r$-neighborhood of point ${\bm{u}}^{(0)}$ 
		\State ${\bm{v}}^{(0)}\ceq\arg\max\lbrace \hat\beta(\bm{z}) ~|~ \bm{z}\in D^{(0)}\rbrace$
		\State ${\bm{w}}^{(0)}\ceq\hat{\bm{\gamma}}({\bm{v}}^{(0)})$
		\State \textbf{assert} $\exists i\in\mathcal{I}: {\bm{w}}^{(0)},{\bm{u}}^{(0)}\in H_i\cap \Gamma(\hat M)$ \Comment If this fails, reduce $r$		
		\While {$\langle \bm{c},{\bm{w}}^{(k)}-{\bm{u}}^{(k)}\rangle > \epsilon_f$}
		\State \textbf{assert} $\exists i\in\mathcal{I}: {\bm{w}}^{(k)},{\bm{u}}^{(k)}\in H_i\cap \Gamma(\hat M)$ \Comment If this fails, reduce $r$
		\State ${\bm{d}}^{(k)}\ceq {\bm{w}}^{(k)}-{\bm{u}}^{(k)}$
		\State $L^{(k)}\ceq\lbrace {\bm{u}}^{(k)} + \lambda {\bm{d}}^{(k)} ~|~ \lambda \in \mathbb{R}_{>0} \rbrace $
		\State ${\bm{u}}^{(k+1)}\ceq \arg\max\lbrace \| \bm{x}-{\bm{u}}^{(k)} \| ~|~ \bm{x}\in L^{(k)} \cap \Gamma(M)\rbrace $
		\State $D^{(k+1)}\ceq H_c({\bm{u}}^{(k+1)}) \cap V_r({\bm{u}}^{(k+1)})$		
		\State ${\bm{v}}^{(k+1)}\ceq\arg\max\lbrace \hat\beta(\bm{z}) ~|~ \bm{z}\in D^{(k+1)} \rbrace$
		\State ${\bm{w}}^{(k+1)}\ceq\hat{\bm{\gamma}}({\bm{v}}^{(k+1)})$
		\State $k\ceq k+1$
		\EndWhile
		\State \textbf{output} ${\bm{u}}^{(k)}$
		\State \textbf{stop}
		\end{algorithmic}\end{algorithm}\end{minipage}
		\begin{multicols}{2}
\noindent The bias $\hat\beta(\bm{z})$ is calculated by equation~\eqref{eq:beta(z)}. The objective projection $\hat{\bm{\gamma}}(\bm{z})$, in equation~\eqref{eq:beta(z)}, is calculated using equations~\eqref{eq:arg_min_beta_i(z)},~\eqref{eq:hat_gamma(z)=gamma_i'(z)} and~\eqref{eq:gamma_i(z)}. 
Step~7 calculates the point ${\bm{w}}^{(0)}$, which is the objective projection of the point ${\bm{v}}^{(0)}$ on the boundary of the recessive polytope. 
Step~8 checks that there exists a recessive half-space $\hat{H}_i$ such that the boundary points ${\bm{w}}^{(0)}$ and ${\bm{u}}^{(0)}$ lie on the hyperplane $H_i$ bounding this half-space. This is necessary so that the movement is carried out on the surface of the recessive polytope, and not through its interior. If this requirement is not met, it is necessary to reduce the radius~$r$ 
of the \mbox{$n$-dimensional} ball~$V_r\left({\bm{u}}^{(0)}\right)$. A suitable $r$ exists by virtue of Proposition~\ref{prp:Gamma(hat_M)}. 
Steps 9--18 implement the main loop of the surface movement method, the geometric interpretation of which is shown in Fig.~\ref{fig:SMM}. This loop is executed while the following condition is true:
\begin{equation}\label{eq:running_criterion}
	\left\langle \bm{c},{\bm{w}}^{(k)} - {\bm{u}}^{(k)}\right\rangle > \epsilon_f,
\end{equation}
where $\epsilon_f$ is a small positive parameter.
Step~10 checks that there exists a recessive half-space $\hat{H}_i$ such that the boundary points ${\bm{w}}^{(k)}$ and ${\bm{u}}^{(k)}$ lie on the hyperplane $H_i$ bounding this half-space. If this requirement is not met, it is necessary to reduce the radius~$r$ 
of the \mbox{$n$-dimensional} ball~$V_r\left({\bm{u}}^{(k)}\right)$.
Step~11 calculates the vector~${\bm{d}}^{(k)}$, which determines the movement direction.

\medskip\noindent
\begin{minipage}{\linewidth}
	\centering
	\includegraphics[width=7 cm]{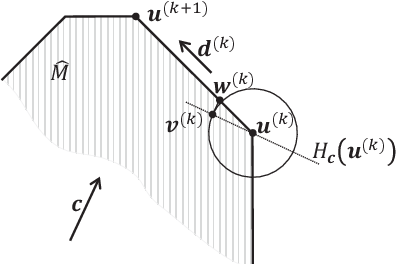}
	\captionof{figure}{Surface movement method.}\label{fig:SMM}
\end{minipage}
\medskip

\noindent Step~12 builds the ray $L^{(k)}$ parallel to the vector ${\bm{d}}^{(k)}$ with the initial point ${\bm{u}}^{(k)}$. 
Step~13 determines the next approximation~${\bm{u}}^{(k+1)}$ as the point belonging to the ray~$L^{(k)}$, which is lying on the boundary of the recessive polytope~$M$ as far as possible from the point~${\bm{u}}^{(k)}$.
Step~14 builds the hyperdisk $D^{(k+1)}$ of radius~$r$ with the center at the point~${\bm{u}}^{(k+1)}$, lying on the hyperplane $H_c\left({\bm{u}}^{(k+1)}\right)$.
Step~15 finds the point ${\bm{v}}^{(k+1)}$ on the hyperdisk~$D^{(k+1)}$ that has the maximum bias. 
Step~16 calculates the point ${\bm{w}}^{(k+1)}$, which is the objective  projection of the point ${\bm{v}}^{(k+1)}$ onto the boundary of the recessive polytope~$\hat{M}$. 
In Step~17, the iteration counter~$k$ is incremented by~1. 
Step~18 passes the control to the beginning of the \texttt{while} loop.
Step~19 outputs the last approximation~ ${\bm{u}}^{(k)}$ as a result. 
Step~20 terminates the algorithm. 

Note that by construction of Algorithm~\ref{alg:SMM}, for any~$k$ the following condition holds:
\begin{equation}\label{eq:u^(k)_in_M_cap_Gamma_hat(M)}
	{\bm{u}}^{(k)}\in M \cap\Gamma(\hat{M}),
\end{equation}
i.e., all points of the sequence~$\left\lbrace {\bm{u}}^{(k)} \right\rbrace$ generated by Algorithm~\ref{alg:SMM} simultaneously lie on the boundary of the feasible polytope $M$ and on the boundary of the recessive polytope $\hat{M}$. Also, it follows from~\eqref{eq:running_criterion} that
\begin{equation}\label{eq:<c,u^(k)><<c,w^k>}
	\left\langle \bm{c},{\bm{u}}^{(k)}\right\rangle 
	< \left\langle \bm{c},{\bm{w}}^{(k)} \right\rangle.
\end{equation}
In addition, the following inequality holds:
\begin{equation}\label{eq:<c,w^k><=<c,u^(k+1)>}
	\left\langle \bm{c},{\bm{w}}^{(k)} \right\rangle
	\leqslant \left\langle \bm{c},{\bm{u}}^{(k+1)} \right\rangle.
\end{equation}

The following lemma guarantees that Algorithm~\ref{alg:SMM} stops in a finite number of iterations.

\begin{lemma}\label{lem:algorithm_1_finiteness}
	Let the feasible polytope~$M$ of the problem~\eqref{eq:LP_problem} be a bounded nonempty set. Then the sequence of points $\left\lbrace {\bm{u}}^{(k)} \right\rbrace$ generated by Algorithm~\ref{alg:SMM} is finite for any $\epsilon_f > 0$.
\end{lemma}
\begin{proof}
		Assume the opposite, that is, Algorithm~\ref{alg:SMM} generates the infinite sequence of points $\left\lbrace {\bm{u}}^{(k)} \right\rbrace$. But then, by virtue of~\eqref{eq:<c,u^(k)><<c,w^k>} and~~\eqref{eq:<c,w^k><=<c,u^(k+1)>}, we obtain the infinite strictly monotonically increasing numerical sequence ${\left\lbrace \left\langle \bm{c},{\bm{u}}^{(k)}\right\rangle \right\rbrace}_{k = 0}^\infty$:
	\begin{equation}\label{eq:monotonically_increasing_sequence}
		\left\langle \bm{c},{\bm{u}}^{(k)}\right\rangle < \left\langle \bm{c},{\bm{u}}^{(k+1)}\right\rangle.  
	\end{equation}	
	Since, by the lemma condition, the feasible polytope~$M$ is a nonempty bounded set, there is a solution $\bar{\bm{x}}$ to the problem LP~\eqref{eq:LP_problem}. By virtue of~\eqref{eq:u^(k)_in_M_cap_Gamma_hat(M)}, the following inequality holds:
	\begin{equation*}
		\left\langle \bm{c},{\bm{u}}^{(k)}\right\rangle \leqslant \left\langle \bm{c},\bar{\bm{x}} \right\rangle
	\end{equation*}
	for all $k=0,1,2,\ldots$ This means that the sequence ${\left\lbrace \left\langle \bm{c},{\bm{u}}^{(k)}\right\rangle \right\rbrace}_{k = 0}^\infty$ is bounded from above. According to the monotone convergence theorem, a monotonically increasing bounded from above sequence converges to its supremum. I.e., there exists $k'\in\mathbb{N}$ such that
	\begin{equation*}
		\forall k>k': \left\langle \bm{c},{\bm{u}}^{(k+1)}\right\rangle - \left\langle \bm{c},{\bm{u}}^{(k)}\right\rangle < \epsilon_f.
	\end{equation*}
	By~\eqref{eq:<c,w^k><=<c,u^(k+1)>} it follows that
	\begin{equation*}
		\forall k>k': \left\langle \bm{c},{\bm{w}}^{(k)}\right\rangle - \left\langle \bm{c},{\bm{u}}^{(k)}\right\rangle < \epsilon_f.
	\end{equation*}
	This is equivalent to
	\begin{equation*}
		\forall k>k': \left\langle \bm{c},{\bm{w}}^{(k)}-{\bm{u}}^{(k)}\right\rangle < \epsilon_f.
	\end{equation*}
	We obtain a contradiction with the condition of the \texttt{while} loop (Step~9 of Algorithm~\ref{alg:SMM}).
\end{proof}

The following theorem shows that, for a sufficiently small~$\epsilon_f$, the sequence of points~$\left\lbrace {\bm{u}}^{(k)} \right\rbrace$ calculated by Algorithm~\ref{alg:SMM} converges to the solution~$\bar{\bm{x}}$ of LP~Problem~\eqref{eq:LP_problem}.

\begin{theorem}\label{trm:convergence}
	Let the feasible polytope~$M$ of the problem~\eqref{eq:LP_problem} be a bounded nonempty set. Let~$\bar{\bm{x}}$ be a solution of LP~Problem~\eqref{eq:LP_problem}. Let~$\left\lbrace \epsilon_\eta \right\rbrace_{\eta = 1}^\infty$ be a monotonically decreasing sequence of positive numbers converging to zero:
	\begin{equation}\label{eq:lim_e_eta=0}
		\mathop {\lim }\limits_{\eta \to \infty } \epsilon_\eta = 0.
	\end{equation}
	Denote by~$u^{(K_\eta)}$ the final point generated by Algorithm~\ref{alg:SMM} with~$\epsilon_f=\epsilon_\eta$ (it exists by virtue of Lemma~\ref{lem:algorithm_1_finiteness}). Then, there exists~$\bar{\eta} \in \mathbb{N}$ such that, for all~$\eta \geqslant \bar{\eta}$, the following equation holds:
	\begin{equation*}
		\left\langle \bm{c}, {\bm{u}}^{(K_\eta)} \right\rangle = \left\langle \bm{c}, \bar{\bm{x}} \right\rangle.
	\end{equation*}
\end{theorem}
\begin{proof}
		Let us show that the sequence $\left\lbrace \left\langle \bm{c}, {\bm{u}}^{(K_{\eta})} \right\rangle\right\rbrace_{\eta = 1}^\infty$ is monotonically increasing. Indeed, it follows from~\eqref{eq:<c,u^(k)><<c,w^k>} and~\eqref{eq:<c,w^k><=<c,u^(k+1)>} that
	\begin{equation}\label{eq:<c,u^k'><=<c,u^k'>}
		\forall k' \leqslant k'':\left\langle \bm{c},{\bm{u}}^{(k')}\right\rangle \leqslant \left\langle \bm{c},{\bm{u}}^{(k'')}\right\rangle.  
	\end{equation}
	By the condition of the theorem,
	\begin{equation*}
		\epsilon_\eta \geqslant \epsilon_{\eta+1}.
	\end{equation*}
	Therefore, by construction of Algorithm~\ref{alg:SMM},
	\begin{equation*}
		K_\eta \leqslant K_{\eta+1}.
	\end{equation*}
	Comparing this with~\eqref{eq:<c,u^k'><=<c,u^k'>}, we obtain
	\begin{equation*}
		\left\langle \bm{c}, {\bm{u}}^{(K_\eta)} \right\rangle \leqslant \left\langle \bm{c}, {\bm{u}}^{(K_{\eta+1})} \right\rangle,
	\end{equation*}
	i.e., the sequence $\left\lbrace \left\langle \bm{c}, {\bm{u}}^{(K_\eta)} \right\rangle\right\rbrace _{\eta = 1}^\infty$ is monotonically increasing.
	Obviously, this sequence is bounded from above by the value of  $\left\langle \bm{c}, \bar{\bm{x}} \right\rangle$. Hence, it has a finite limit:
	\begin{equation*}
		\mathop {\lim }\limits_{\eta \to \infty } \left\langle \bm{c}, {\bm{u}}^{(K_\eta)} \right\rangle = \bar{f}.
	\end{equation*}
	Algorithm~\ref{alg:SMM}, within each iteration, passes\footnote{The passage of a polytope face/edge is understood as movement inside a linear manifold of dimension $k$ in the presence of $k$ degrees of freedom.} one face/edge of the recessive polytope~$\hat{M}$ in the direction of maximizing the value of the objective function. In this case, each face/edge is traversed no more than once, since the polytope~$\hat{M}$ is a convex set. This means that there exists $\bar\eta \in\mathbb{N}$ such that for all $\eta \geqslant \bar\eta$ the following equation holds:
	\begin{equation*}
		{\bm{u}}^{(K_\eta)} = {\bm{u}}^{\left(K_{\bar\eta}\right)},
	\end{equation*}
	and
	\begin{equation*}
		\left\langle \bm{c}, {\bm{u}}^{(K_\eta)} \right\rangle = \bar{f}.
	\end{equation*}
	By construction of Algorithm~\ref{alg:SMM}, taking into account~\eqref{eq:lim_e_eta=0}, this is possible only when
	\begin{equation}\label{eq:c,w-u^bar_K=0}
		\left\langle \bm{c}, {\bm{w}}^{\left(K_{\bar{\eta}}\right)} - {\bm{u}}^{\left(K_{\bar{\eta}}\right)} \right\rangle = 0.
	\end{equation}
	Let us show that, in this case,
	\begin{equation*}
		\left\langle \bm{c}, {\bm{u}}^{\left(K_{\bar{\eta}}\right)} \right\rangle = \left\langle \bm{c}, \bar{\bm{x}} \right\rangle,
	\end{equation*}
	i.e., the point ${\bm{u}}^{\left(K_{\bar{\eta}}\right)}$ is a solution of  LP~Problem~\eqref{eq:LP_problem}. Denote $\bm{u}'={\bm{u}}^{\left(K_{\bar{\eta}}\right)}$, and assume the opposite, namely that there exists a point
	\begin{equation}\label{eq:u''_in_M} 
		\bm{u}'' \in M
	\end{equation}
	such that
	\begin{equation*}
		\left\langle \bm{c},{\bm{u}}''\right\rangle >\left\langle c, {\bm{u}}'\right\rangle.
	\end{equation*}
	This is equivalent to
	\begin{equation}\label{eq:<c,u''-<c,u'>>0}
		\left\langle \bm{c},{\bm{u}}''-{\bm{u}}'\right\rangle > 0.
	\end{equation}
	Fig.~\ref{fig:Theorem1} illustrates the following part of the proof.

\medskip\noindent
\begin{minipage}{\linewidth}
	\centering
	\includegraphics[width=7 cm]{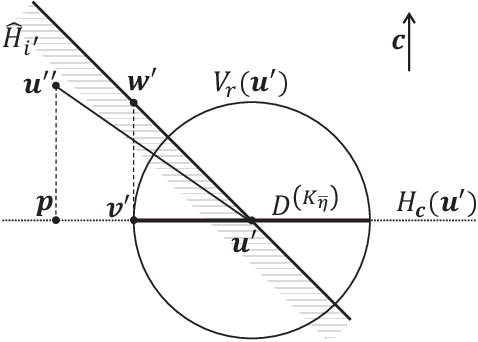}
	\captionof{figure}{Illustration to proof of Theorem \ref{trm:convergence}.}\label{fig:Theorem1}
\end{minipage}
\medskip

	Based on Definition~\ref{def:H_c}, we can calculate the orthogonal projection~${\bm{p}}$ of the point~${\bm{u}}''$ onto the objective hyperplane~$H_c\left({\bm{u}}'\right)$ passing through the point~${\bm{u}}'$:
	\begin{equation}\label{eq:p}
		{\bm{p}}={\bm{u}}'' - \frac{\left\langle \bm{c},{\bm{u}}''-{\bm{u}}'\right\rangle}{{\left\|\bm{c}\right\|}^2}\bm{c}.
	\end{equation}
	Note that
	\begin{equation}\label{eq:||p-u'||>0}
		\left\| {\bm{p}}-{\bm{u}}'\right\|\ne 0,
	\end{equation}
	since otherwise, according to Definition~\ref{def:recessive_halfspace}, the point~${\bm{u}}''$ cannot belong to the recessive polytope~$\hat{M}$, which contradicts~\eqref{eq:u''_in_M}. Choose~$r\in\mathbb{R}$ satisfying the condition
	\begin{equation}\label{eq:r>0}
		r>0,
	\end{equation}
	for which there is~$i'\in\mathcal{I}$ such that
	\begin{equation}\label{eq:hat_gamma(v),u'_in_H_i'_cap_Gamma(M)}
		\hat{\bm{\gamma}}({\bm{v}'}), {\bm{u}}' \in H_{i'} \cap \Gamma(M),
	\end{equation}
	where
	\begin{equation}\label{eq:v}
		{\bm{v}'}={\bm{u}}'+\frac{r}{\left\| {\bm{p}}-{\bm{u}}'\right\| }({\bm{p}}-{\bm{u}}').
	\end{equation}
	This is possible by virtue of Proposition~\ref{prp:Gamma(M)}. Without loss of generality, we can assume that~$r$ satisfies all asserts in the case when $\epsilon_f = \epsilon_{\bar{\eta}}$. Thus
	\begin{equation}\label{eq:v'_in_D_K_bar_eta}
		\bm{v}' \in D^{(K_{\bar\eta})}, 
	\end{equation}
	where $D^{K_{\bar\eta}}$ is the disk constructed in Step~14 of algorithm~\ref{alg:SMM} with $k=K_{\bar\eta}-1$. Note that
	\begin{equation*}
		i' = \arg\min \left\{ \beta_i(\bm{v}')   \left| i\in\mathcal{I} \right.\right\},	
	\end{equation*}
	since otherwise $\hat{\bm{\gamma}}({\bm{v}'}) \notin H_{i'}$, which contradicts~\eqref{eq:hat_gamma(v),u'_in_H_i'_cap_Gamma(M)}. According to Proposition~\ref{prp:hat_gamma(z)=gamma_i'(z)}, it follows that
	\begin{equation}\label{eq:hat_gamma(v)=gamma_i'(v)}
		\hat{\bm{\gamma}}({\bm{v}'})=\bm{\gamma}_{i'}({\bm{v}'}).
	\end{equation}
	Set
	\begin{equation*}
		\bm{w}'=\hat{\bm{\gamma}}({\bm{v}'}).
	\end{equation*}
	Taking into account~\eqref{eq:hat_gamma(v)=gamma_i'(v)}, the following equation holds:
	\begin{equation*}
		\bm{w}'=\bm{\gamma}_{i'}({\bm{v}'}).
	\end{equation*}
	Using Proposition~\ref{prp:c-projection}, we obtain
	\begin{equation}\label{eq:w}
		{\bm{w}'}={\bm{v}'} - \frac{\left\langle \bm{a}_{i'},{\bm{v}'}\right\rangle - b_{i'}}{\left\langle \bm{a}_{i'},\bm{c}\right\rangle}\bm{c}.
	\end{equation}
	Since ${\bm{u}}'\in H_{i'}$, it follows from~\eqref{eq:H_i} that
	\begin{equation}\label{eq:<a_i',u'>=b_i'}
		\left\langle  \bm{a}_{i'},{\bm{u}}' \right\rangle = b_{i'}.
	\end{equation}
	Therefore,~\eqref{eq:w} can be rewritten in the form
	\begin{equation*}
		{\bm{w}'}={\bm{v}'} - \frac{\left\langle \bm{a}_{i'},{\bm{v}'}\right\rangle - \left\langle  \bm{a}_{i'},{\bm{u}}' \right\rangle}{\left\langle \bm{a}_{i'},\bm{c}\right\rangle}\bm{c}.
	\end{equation*}
	Substituting the right side of equation~\eqref{eq:v} instead of~${\bm{v}'}$, we obtain from here
	\begin{multline*}
	{\bm{w}'}={\bm{u}}'+\frac{r}{\left\| {\bm{p}}-{\bm{u}}'\right\| }({\bm{p}}-{\bm{u}}') - \\
	- \frac{\left\langle \bm{a}_{i'},{\bm{u}}'+\frac{r}{\left\| {\bm{p}}-{\bm{u}}'\right\| }({\bm{p}}-{\bm{u}}')\right\rangle - \left\langle  \bm{a}_{i'},{\bm{u}}' \right\rangle}{\left\langle \bm{a}_{i'},\bm{c}\right\rangle}\bm{c}.
\end{multline*}
	which is equivalent to
	\begin{equation}\begin{gathered} \label{eq:w2}
		{\bm{w}'}={\bm{u}}'+\frac{r}{\left\| {\bm{p}}-{\bm{u}}'\right\| }({\bm{p}}-{\bm{u}}') - \hfill \\ 
		-\frac{\left\langle \bm{a}_{i'},\frac{r}{\left\| {\bm{p}}-{\bm{u}}'\right\| }({\bm{p}}-{\bm{u}}')\right\rangle}{\left\langle \bm{a}_{i'},\bm{c}\right\rangle}\bm{c}.
	\end{gathered}\end{equation} 
	According to~\eqref{eq:half-space} we have
	\begin{equation}\label{eq:hat_H_i'}
		\hat H_{i'}=\left\lbrace \bm{x}\in\mathbb{R}^n \middle|\left\langle  \bm{a}_{i'},\bm{x} \right\rangle \leqslant b_{i'} \right\rbrace.
	\end{equation}
	Using~\eqref{eq:<a_i',u'>=b_i'}, equation~\eqref{eq:hat_H_i'} can be rewritten in the form
	\begin{equation}\label{eq:b_i'=<a_i',u'>}
		\hat H_{i'}=\left\lbrace \bm{x}\in\mathbb{R}^n \middle|\left\langle  \bm{a}_{i'},\bm{x} \right\rangle \leqslant \left\langle  \bm{a}_{i'},{\bm{u}}' \right\rangle \right\rbrace.
	\end{equation}
	From~\eqref{eq:u''_in_M}, it follows that ${\bm{u}}''\in \hat H_{i'}$. Comparing this with~\eqref{eq:b_i'=<a_i',u'>}, we obtain
	\begin{equation*}
		\left\langle  \bm{a}_{i'},{\bm{u}}'' \right\rangle \leqslant \left\langle  \bm{a}_{i'},{\bm{u}}' \right\rangle,
	\end{equation*}
	which is equivalent to
	\begin{equation}\label{eq:<a_i',u'-u''>>=0}
		\left\langle  \bm{a}_{i'},{\bm{u}}'-{\bm{u}}'' \right\rangle \geqslant 0.
	\end{equation}
	
	Since the half-space $\hat{H}_{i'}$ is recessive, according to Proposition~1 in~\cite{Sokolinsky2023:TR}, the following inequality holds
	\begin{equation}\label{eq:<a_i',c>>0}
		\left\langle \bm{a}_{i'}, \bm{c} \right\rangle > 0.
	\end{equation}
	
	By virtue of~\eqref{eq:w2} and~\eqref{eq:p} we have
	\begin{equation*}
		\begin{array}{l}
			\left\langle \bm{c},\bm{w}'-{\bm{u}}'\right\rangle = r\frac{\left\langle \bm{c},{\bm{p}}-{\bm{u}}'\right\rangle - \frac{\left\langle \bm{a}_{i'},{\bm{p}}-{\bm{u}}'\right\rangle}{\left\langle \bm{a}_{i'},\bm{c}\right\rangle} \left\| \bm{c} \right\|  }{\left\| {\bm{p}}-{\bm{u}}'\right\| } \\
			=  \frac{r\left\| \bm{c} \right\|}{\left\| {\bm{p}}-{\bm{u}}'\right\|\left\langle \bm{a}_{i'},\bm{c}\right\rangle }\left( - \left\langle \bm{a}_{i'}, {\bm{u}}'' - \frac{\left\langle \bm{c},{\bm{u}}''-{\bm{u}}'\right\rangle}{{\left\|\bm{c}\right\|}^2} \bm{c} - {\bm{u}}'\right\rangle \right) \\
			=  \frac{r\left\| \bm{c} \right\|}{\left\| {\bm{p}}-{\bm{u}}'\right\|\left\langle \bm{a}_{i'},\bm{c}\right\rangle }\left( \left\langle \bm{a}_{i'}, {\bm{u}}' - {\bm{u}}'' \right\rangle + \frac{\left\langle \bm{c},{\bm{u}}''-{\bm{u}}'\right\rangle\left\langle \bm{a}_{i'}, \bm{c}\right\rangle}{{\left\|\bm{c}\right\|}^2} \right)
		\end{array} 
	\end{equation*}
	According to \eqref{eq:r>0}, \eqref{eq:||p-u'||>0}, \eqref{eq:<a_i',c>>0}, \eqref{eq:<a_i',u'-u''>>=0} and \eqref{eq:<c,u''-<c,u'>>0}, it follows that
	\begin{equation*}
		\left\langle \bm{c},{\bm{w}'}-{\bm{u}}'\right\rangle > 0.
	\end{equation*}
	Recalling that ${\bm{u}}'={\bm{u}}^{\left(K_{\bar{\eta}}\right)}$, we rewrite the last inequality in the form
	\begin{equation}\label{eq:<c,w'-u_K_bar_eta>0}
		\left\langle \bm{c},{\bm{w}'} - {\bm{u}}^{\left(K_{\bar{\eta}}\right)}\right\rangle > 0.
	\end{equation}
	Taking in account Proposition~\ref{prp:metrics}, by construction of Algorithm~\ref{alg:SMM} (steps~15,~16), we have
	\begin{equation*}
		\left\langle \bm{c},{\bm{w}}^{\left(K_{\bar{\eta}}\right)}\right\rangle \geqslant \left\langle \bm{c},{\bm{w}'}\right\rangle.
	\end{equation*}
	Comparing this with~\eqref{eq:<c,w'-u_K_bar_eta>0}, we obtain
	\begin{equation*}
		\left\langle \bm{c},{\bm{w}}^{\left(K_{\bar{\eta}}\right)} - {\bm{u}}^{\left(K_{\bar{\eta}}\right)}\right\rangle > 0.
	\end{equation*}
	We have thus reached a contradiction with~\eqref{eq:c,w-u^bar_K=0}.
\end{proof}

\section{Discussion}\label{sec:discussion}

This section discusses the strengths and weaknesses of the proposed method, as well as reveals the ways of its practical implementation based on the synthesis of high-performance computing and neural network technologies.

First of all, let us look at the issue concerning the implementation of the algorithm~\ref{alg:SMM} in the form of a computer program. Steps~15 of Algorithm~\ref{alg:SMM} should calculate the point of hyperdisk~$D^{(k+1)}$ having the maximum bias. This task can be solved using the approach described in article~\cite{Sokolinsky2023:TR}. It is based on the pseudoprojection operation, which is a generalization of the metric projection to a convex closed set. The idea of the solution is as follows. For the current approximation $\bm{u}$, let us find the numbers of  all hyperplanes passing through $\bm{u}$:
\begin{equation*}
	\forall i \in \mathcal{H} : \left\langle  \bm{a}_i,\bm{u} \right\rangle = b_i;
\end{equation*}
\begin{equation*}
	\forall i \in \mathcal{P}\backslash\mathcal{H} : \left\langle  \bm{a}_i,\bm{u} \right\rangle \ne b_i.
\end{equation*}
Let $\mathfrak{H}$ be the set of all subsets of the set $\mathcal{H}$, except the empty set. This set defines a set of linear manifolds of the form
\begin{equation*}
	V(\mathcal{S}) = \bigcap\limits_{i \in \mathcal{S}} H_i, 		
\end{equation*}	
where $\mathcal{S} \in \mathfrak{H}$.
Then, any face $Q$ of the polytope $M$ containing $\bm{u}$ can be represented as follows: 
\begin{equation*}
	Q = M \cap V(\mathcal{S})
\end{equation*}
with $\mathcal{S} \in \mathfrak{H}$.
Given the gradient $\bm{c}$ of the objective function, we can find the unit direction vector of maximum increase of the objective function for each face using the pseudoprojection operation. The face with the maximum increment will provide us with the direction vector $d$ that will allow us to find the next approximation (see Fig.~\ref{fig:SMM}). We plan to describe this method in detail in a separate paper. However, we have already performed and tested a parallel implementation of this method. The source codes and results of the runs are freely available on GitHub at \url{https://github.com/leonid-sokolinsky/BSF-Surface-movement-method}. The drawback of this method is that its time complexity can be estimated as $O(2^m)$, where $m$ is the number of constraints of the LP problem. We see the following way to solve this issue involving an artificial neural network. Using the approach described in paper~\cite{Olkhovsky2022:TR}, we replace the hyperdisk $D^{(k+1)}$ in Algorithm~\ref{alg:SMM} with a set of points called the receptive field. We map each point of the receptive field to its bias relative to the boundary of the feasible polytope. As a result, we obtain a matrix of dimension $(n-1)$, which is a local image of the LP problem. The locality of the image means that we obtain a visual representation of the surface not of the entire feasible polytope, but only of some part of it in the neighborhood of the point of the current approximation. This image is fed to the input of a pre-trained feed forward neural network, which outputs the vector~$\bm{d}$ indicating the direction of movement on the surface of the feasible polytope towards the maximum increase in the value of the objective function. Denote by~$\mathfrak{G}({\bm{u}})$ the function that constructs a receptive field centered at the point~${\bm{u}}$ and calculates the local image of the LP problem at this point. The algorithm for constructing a multidimensional image of the LP problem is described and investigated in~\cite{Olkhovsky2022:TR}. Denote by~$\operatorname{DNN}$ a deep neural network, which consumes a local image of the LP problem and outputs the vector~$\bm{d}$ that defines the direction of movement along the surface of the feasible polytope. Thus, Algorithm~\ref{alg:SMM} can be transformed into Algorithm~\ref{alg:DNN}. The set of labeled examples required for DNN training can be obtained using the method described above. We plan to conduct a separate study specifically addressing this subject.

\end{multicols}
\begin{minipage}{0.95\textwidth}\begin{algorithm}[H]\caption{LP method using DNN}\label{alg:DNN}\begin{algorithmic}[1]
		\Require \begin{flushleft} $\hat H_i=\left\lbrace \bm{x}\in\mathbb{R}^n \middle|\left\langle  \bm{a}_i,\bm{x} \right\rangle \leqslant b_i \right\rbrace;\; \hat M=\bigcap\limits_{i \in \mathcal{I}} {\hat H_i}$ \end{flushleft}
		\State \textbf{input} ${\bm{u}}^{(0)}$
		\State \textbf{assert} ${\bm{u}}^{(0)}\in M \cap \Gamma(\hat{M})$
		\State $k\ceq 0$
		\State $D^{(0)}\ceq \mathfrak{G}({\bm{u}}^{(0)})$
		\State ${\bm{d}}^{(0)}\ceq \operatorname{DNN}(D^{(0)})$
		\While {${\bm{d}}^{(k)} \ne \mathbf{0}$}
		\State $L^{(k)}=\lbrace {\bm{u}}^{(k)} + \lambda {\bm{d}}^{(k)} ~|~ \lambda \in \mathbb{R}_{>0} \rbrace $
		\State ${\bm{u}}^{(k+1)}\ceq \arg\max\lbrace \| \bm{x}-{\bm{u}}^{(k)} \| ~|~ \bm{x}\in L^{(k)} \cap \Gamma(M)\rbrace $
		\State $D^{(k+1)}\ceq \mathfrak{G}({\bm{u}}^{(k+1)})$		
		\State ${\bm{d}}^{(k+1)}\ceq \operatorname{DNN}(D^{(k+1)})$
		\State $k\ceq k+1$
		\EndWhile
		\State \textbf{output} ${\bm{u}}^{(k)}$
		\State \textbf{stop}
\end{algorithmic}\end{algorithm}\end{minipage}\begin{multicols}{2}

Let us estimate the time complexity of Algorithm~\ref{alg:DNN}. The surface movement method visits\footnote{A visit to a hyperplane is understood as a rectilinear movement from the point of entry to the hyperplane to the point of the first change in the direction of movement.} each hyperplane of the recessive polytope no more than once. A visit to a hyperplane is performed within one iteration of the \texttt{while} loop (steps 6--12 of Algorithm~\ref{alg:DNN}). Therefore, the total number of iterations can be estimated as $O(m)$, where~$m$ is the number of constraints of LP~Problem~\eqref{eq:LP_problem}. Finding the next approximation~${\bm{u}}^{(k+1)}$ in Step~8 of Algorithm~\ref{alg:DNN} can be implemented by dichotomy. Thus, the number of operations in Step~8 does not depend on the number of constraints~$m$, nor on the dimension~$n$, and, for large values of~$m$ and~$n$, it can be estimated as a constant. The most time-consuming operation is the construction of a local image of LP~Problem~\eqref{eq:LP_problem} in Step~9 of Algorithm~\ref{alg:DNN}. A~receptive field in the form of a hypercubic lattice consists of~$\eta^{(n-1)}$ points, where~$n$ is the dimension of space, and~$\eta$ is the number of points in one dimension. However, recent study~\cite{Olkhovsky2023:Rus} has shown that a cruciform receptive field with the number of points~$(\eta-1)n+1$ gives results that are not inferior to the hypercubic one in terms of the accuracy of solving the LP problem. A pre-trained feedforward neural network DNN calculates the movement direction vector $\bm{d}$ at step~10 in a time that depends only on~$n$, since it is fed $(\eta-1)n+1$ numbers (in the case of a cruciform receptive field). Thus, the total time complexity of the algorithm can be estimated as $O(mn)$.

In addition, we note that the scalability boundary\footnote{The scalability boundary refers to the number of processor nodes of a cluster computing system on which the maximum speedup is achieved.} of the algorithm for constructing a visual image of the LP problem can be estimated as $O(\sqrt{2n^2m+m^2n+8nm-6m})$~\cite{Olkhovsky2022:TR}. Under the assumption that $m=O(n)$, we obtain, for the scalability boundary, the estimation $O(n\sqrt{n})$, which is close to a linear relationship. This means that the algorithm for constructing a local image of the LP problem can be effectively parallelized on a large number of processor nodes of a cluster computing system. So for $n=7$ and $m=15$, computational experiments show a peak of speedup on 326 processor nodes~\cite{Olkhovsky2022:TR}. Note that the number of iterations of Algorithm~\ref{alg:DNN} does not depend on $\epsilon_f$, since there is no such parameter in this algorithm.

The surface movement method is self-correcting. Therefore, it can potentially be used to solve non-stationary problems. At that, if only the objective function changes, then Algorithm~\ref{alg:DNN} does not require any crucial changes at all. It is important that the rate of correction outperforms the rate of change. If the system of constraints changes (without changing the dimension), then Algorithm~\ref{alg:DNN} will require certain modifications, since the current approximation may ``dive'' into the polytope or ``break away'' from its surface. The authors intend to study this issue in detail in the future.

\begin{figure*}[!ht]
	\centering
	\begin{subfigure}[t]{\columnwidth}
		\centering
		\includegraphics[width=5 cm]{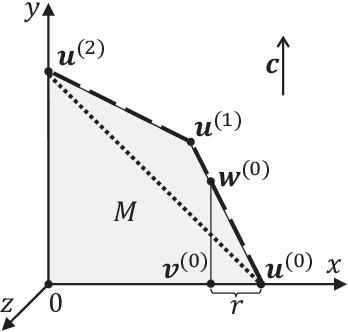}
		\captionof{figure}{Optimal objective path is indicated by dashed line;
			path of minimum length is shown by dotted line.}\label{fig:Optimal-objective-path}
	\end{subfigure}
	\hfill
	\begin{subfigure}[t]{\columnwidth}
		\centering
		\includegraphics[width=5 cm]{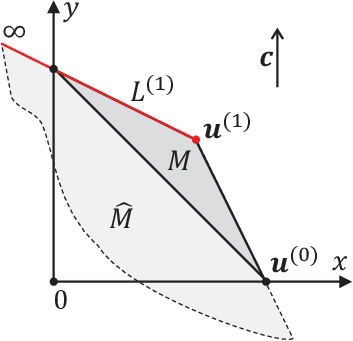}
		\captionof{figure}{Optimal objective path goes from ${\bm{u}}^{(1)}$ to infinity on recessive polytope $\hat{M}$.}\label{fig:To-Infinity}
	\end{subfigure}
	\caption{Special cases of optimal objective path.}
\end{figure*}
\medskip

Algorithm~\ref{alg:DNN} can also be used to solve LP problems in real time. Indeed, the number of iterations is bounded by the parameter~$m$. With a fixed $n$, building a local image of the LP problem requires a fixed number of operations. At the same time, the image construction procedure is effectively parallelized on a large number of processor nodes. The feedforward neural network DNN processes the local image of the LP problem in a fixed time, depending only on~$n$. The work of a neural network can also be efficiently parallelized using GPUs. In the future, it is possible to refuse the travelling along the faces/edges of the feasible polytope and analyze with the use of a neural network the image of the entire feasible polytope obtained from the apex point (see~\cite{Sokolinsky2023:TR}). The neural network will produce an approximate solution, which can be refined by analyzing a fixed number of local images with increasing detail. However, this issue needs further research.

The Algorithm~\ref{alg:SMM} constructs the optimal objective path to a solution of the LP problem on the surface of the feasible polytope. This directly follows from the construction of the algorithm and Proposition~\ref{prp:metrics}. An interesting question is whether this path is always the path of the shortest length in the sense of the Euclidean metric. The following simple example in the space $\mathbb{R}^3$ shows that this is not always the case. For the system of constraints
\begin{equation*}
	\left\lbrace 
	\begin{array}{l}
		x+2y \leqslant 2,\\
		2x+y \leqslant 2,\\
		x \geqslant 0,\\
		y \geqslant 0,\\
		z = 0;
	\end{array} 
	\right. 
\end{equation*}
and the objective function $f(x,y,z)=y$, the optimal objective path may not coincide with the shortest path (see Fig.~\ref{fig:Optimal-objective-path}).

The question may also arise, whether it is possible to replace $\Gamma(M)$ with $\Gamma(\hat{M})$ in Step~13 of Algorithm~\ref{alg:SMM}, since this reduces the number of inequalities used for checking the condition $x\in\Gamma(M)$. The answer is negative. For example, in the space $\mathbb{R}^2$, for the system of constraints
\begin{equation*}
	\left\lbrace 
	\begin{array}{l}
		x+2y \leqslant 2,\\
		2x+y \leqslant 2,\\
		x+y \geqslant 1,\\
		x \geqslant 0,\\
		y \geqslant 0;\\
	\end{array} 
	\right. 
\end{equation*}
and the objective function $f(x,y)=y$, with ${\bm{u}}^{(1)} = \left( \tfrac{2}{3}, \tfrac{2}{3} \right)$, we obtain \[\max\lbrace \| \bm{x}-{\bm{u}}^{(1)} \| ~|~ \bm{x}\in L^{(1)} \cap \Gamma(\hat{M})\rbrace = +\infty\] (see Fig.~\ref{fig:To-Infinity}).

As a drawback of the surface movement method, it can be noted that it is not affine invariant, since the cost functional is identified with the vector. Thus, the behavior of the method depends on the Euclidean structure defined by the coordinate system.

The scientific contribution and theoretical significance of the proposed method is that, for the first time, it opens up the possibility of using feedforward neural networks to solve multidimensional LP problems based on the analysis of their images.

\section{Conclusion}\label{sec:conclusion}

The article presented a new method for solving the linear programming problem (LP), called the ``surface movement method''. This method constructs, on the surface of the polytope bounding the feasible region, a path from an initial point to a point of solving the LP problem. The movement vector is always constructed in the direction of the maximum increase/decrease in the value of the objective function. The resulting path is called the optimal objective path.

The surface movement method assumes the use of a feedforward neural network to determine the direction of movement along the faces of the feasible polytope. To do this, a local image of the multidimensional LP problem is constructed at the point of the current approximation, which is fed to the input of the neural network. The set of labeled precedents needed for training a neural network can be obtained using the apex method.

To build a theoretical basis of the surface movement method, the concept of the objective projection is introduced. The objective projection is an oblique projection in the direction parallel to the gradient vector of the objective function. A scalar quantity called bias is defined. The bias modulus is equal to the distance from the point to its objective projection. The bias sign is determined by the position of the point inside or outside of the feasible polytope. An equation for calculating the bias of a point relative to the boundary of the feasible polytope is obtained. It is shown that a larger bias corresponds to a larger value of the objective function. A formalized description of the surface movement method in the form of an algorithm is presented. The main convergence theorem of the surface movement method to the solution of the LP problem in a finite number of iterations is proved. A version of the surface movement algorithm using a function of constructing a local multidimensional image of the LP problem and a deep neural network is provided.

As directions for further research, the following can be indicated.
\begin{enumerate}
	\item Design and training of a DNN network capable of calculating the movement vector in the direction of maximizing the value of the objective function for multidimensional LP problems.
	\item Development of a software package for a cluster computing system implementing Algorithm~\ref{alg:DNN} by combining supercomputer and neural network technologies.
	\item Study of the dependence of the DNN network accuracy on the density of the receptive field.
	\item Study of the suitability of the surface movement method for solving non-stationary LP problems.
	\item Study of the suitability of the surface movement method for solving LP problems in real time.
	\item Development of a new visual method for solving LP problems using neural networks based on the analysis of the image of the feasible polytope as a whole.
\end{enumerate}

\section*{Acknowledgments}
\noindent The research was supported by Russian Science Foundation (project No.~23-21-00356).

\section*{Appendix:  Notations}
\begin{tabbing}
MM. \= M \= MMMMMMMMMMMMM \kill
$\mathbb{R}^n$\> real Euclidean space\\
$\left\|\cdot\right\|$\> Euclidean norm\\
$\langle\cdot,\cdot\rangle$\> dot product of two vectors\\
$f(\bm{x})$\> linear objective function \\
$\bm{c}$\> gradient of objective function $f(\bm{x})$ \\
$\bm{e_c}$\> unit vector parallel to vector $\bm{c}$ \\
$\bar{\bm{x}}$\> solution of LP problem \\
$\bm{a}_i$\> $i$th row of matrix $A$ \\
$H_i$\> hyperplane defined by equation $\left\langle  \bm{a}_i,\bm{x} \right\rangle = b_i$ \\
$\hat H_i$\> half-space defined by inequality $\left\langle  \bm{a}_i,\bm{x} \right\rangle \leqslant b_i$ \\
$\mathcal P$\> set of row indexes of matrix $A$ \\
$\mathcal{I}$\> index set of recessive half-spaces \\
$M$\> feasible polytope \\
$\hat{M}$\> recessive polytope \\
$\Gamma(M)$\>\> boundary of $M$ \\
$\Gamma(\hat{M})$\>\> boundary of $\hat{M}$ \\
$\bm{\gamma}_i(\bm{z})$\>\> objective projection of $\bm{z}$ onto $H_i$\\
$\beta_i(\bm{z})$\>\> bias of $\bm{z}$ relative to $H_i$\\
$\hat{\bm{\gamma}}(\bm{z})$\>\>  objective projection of $\bm{z}$ onto $\Gamma(\hat M)$ \\
$\hat{\beta}(\bm{z})$\>\>  bias of $\bm{z}$ relative to $\Gamma(\hat M)$ \\
$V_r(\bm{x})$\>\> $n$-dimensional ball of radius $r$ centered \\
	\>\>at point $\bm{x}$
\end{tabbing}

\bibliographystyle{spmpsci}
\bibliography{Bibliography}

\vspace{1.5cm}

\small
\noindent
\textit{\textbf{Nikolay~A.~Olkhovsky} is a postgraduate student at the School of Electronics and Computer Science of South Ural State University. His field of research is computational mathematics, parallel computing, and artificial neural networks.}\\
\includegraphics[height=10pt, width=10pt]{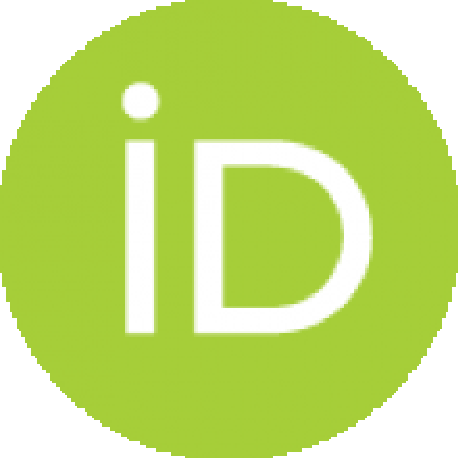} \url{https://orcid.org/0009-0008-9078-4799}

\vspace{1.5cm}

\noindent
\textit{\textbf{Leonid B. Sokolinsky} is the head of the System Programming Department at South Ural State University. He received a degree of candidate of sciences in physics and mathematics from Moscow State University in 1990. In 2003, he received a degree of doctor of sciences in physics and mathematics from Chelyabinsk State University. He has published more than 150 scientific papers. His H-index in Scopus is~10. His fields of research are machine learning, parallel computing, database systems and computational mathematics.}\\
\includegraphics[height=10pt, width=10pt]{orcid_128x128.eps} \url{https://orcid.org/0000-0001-9997-3918}
    
\end{multicols*}
  
\end{document}